\documentclass[12pt]{article}
\usepackage{amssymb,amsmath,amsthm,amscd,latexsym}
\usepackage{mathrsfs}
\usepackage{mathrsfs}
\usepackage{amsfonts}
\usepackage{amsmath}
\usepackage{amssymb}
\usepackage{amscd}
\usepackage{mathrsfs}
\usepackage{amssymb}
\usepackage{amsmath}
\usepackage{amsthm}
\usepackage{color}
\usepackage{latexsym}
\usepackage{indentfirst}
\usepackage{anysize}
\usepackage[normalem]{ulem}

\renewcommand{\paragraph}{\roman{paragraph}}
\input cyracc.def

\newtheorem{theorem}{Theorem}

\newtheorem{lem}{Lemma}

\newtheorem{thm}{Theorem}

\newtheorem{ex}{Example}

\theoremstyle{definition}
\newtheorem{definition}[theorem]{Definition}

\newcommand{\R}{\mathbb{R}}

\begin{document}
%\begin{CJK*}{GBK}{song}\CJKtilde
\title{\bf The covering radius of permutation designs }

\author{
\small{ Patrick Sol\'e$^{1,2}$}\\ 
\small{${}^1$I2M,(Aix-Marseille Univ., Centrale Marseille, CNRS), Marseilles, France}\\
\small{${}^2$ Corresponding Author}
}

%Key Laboratory of Intelligent Computing and Signal Processing of Ministry of Education, School of Mathematics Sciences, Anhui University, Hefei, 230601, China.

\date{}
\maketitle
\begin{abstract} A notion of $t$-designs in the symmetric group on $n$ letters was introduced by Godsil in 1988. In particular $t$-transitive sets of permutations form a $t$-design. 
We derive  upper bounds on the covering radius of these designs, as a function of $n$ and $t$ and in terms of the largest zeros of Charlier polynomials.

\end{abstract}

{\bf Keywords:}  permutation designs, Charlier polynomials\\

{\bf AMS Math Sc. Cl. (2010):} Primary 05E35, Secondary O5E20,  05E24 
%%%%%%%%%%%%%%%%%%%%%%%%%%%%%%%%%%%%%%%%%%%%%%%%%%%%%%%%%
\section{Introduction}
Packing codes in the symmetric group for the Hamming distance have been studied since the 1970s \cite{BCD}. See \cite{C} for a recent survey.
In \cite{CW} {\em covering codes} in the symmetric group are considered. In particular it is shown there that $t$-transitive groups in the symmetric group $S_n$ on $n$ letters have covering radius
for the Hamming distance at most $n-t,$ and that this bound is tight.

In the present paper, we prove bounds of similar order for $t$-designs in $S_n$ in the sense of Godsil \cite{CG,G,Gp}. These objects are defined in the setting of polynomial spaces, a generalization
of association schemes \cite{BI,D}.
An alternative definition in the language of distance degree regular spaces can be found in \cite{SRS,SXS}. It is known that $t$-transitive groups  are $t$-designs, but the converse is not generally true.
To derive these bounds we extend the method of \cite{S} from the Hamming space to the space of permutations. 
This method is based on the polynomials orthogonal w.r.t. the weight distribution of the cosets of the  code considered. If the dual distance of the code is large enough, these polynomials coincide
with the celebrated  Krawtchouk polynomials \cite{Sz}, and the zeros of these can be used to bound the extreme points of the distribution.
In that coding analogy, Charlier polynomials and their zeros play the role of Krawtchouk polynomials.
An important difference between these two families of orthogonal polynomials is that the integrality of the zeros of Charlier polynomials is easier to decide. This technical point simplifies 
the proofs of the bounds in comparison with the coding situation.

The note is arranged as follows. The next section collects notions and definitions needed for the other sections. Section 3 recalls the current results on $t$-designs of permutations.
Section 4 contains the main results.
%%%%%%%%%%%%%%%%%%%%%%%%%%%%%%%%%%%%%%%%%%%%%%%%%%%
\section{Background material}
\subsection{Permutations groups}

A permutation group $G$ acting on a set $X$ of $n$ elements is transitive if there is only one orbit on $X.$ It is $t$-transitive if it is transitive in its action on ${ X \choose t}$ the set of distinct $t$-tuples from $X.$
It is sharply $t$-transitive if this action is regular, concretely if $|G|=\frac{n!}{(n-t)!}.$ We extend this terminology by relaxing the group hypothesis to a set of permutations action on $X.$
It is well-known amongst geometers and group theorists  that a set of sharply $2$-transitive permutations on a set of size $n$
is equivalent to the existence of a projective plane $PG(2,n),$ that is to say a
$2-(n^2+n+1,n+1,1)$ design \cite{C}.

\subsection{Permutation codes }

 Consider the {\em symmetric group} on $n$ letters $S_n$ with metric
$$d_S(\sigma,\theta)=n-F(\sigma \theta^{-1}),$$
where $F(\nu)$ denotes the number of fixed points of $\nu.$
The space $(S_n,d_S)$ is a  metric space.
Let $w_k$ denote the numbers of permutations on $n$ letters with $k$ fixed points. A generating function for these numbers (sometimes called rencontres numbers) is
$$\sum_{k=0}^nw_ku^k=n!\sum_{j=0}^n\frac{(u-1)^j}{j!},$$ as per \cite{w}.  It is clear that $d_S$ is not a shortest path distance since $ d_S(\sigma,\theta)=1$ is impossible.
Codes in $(S_n,d_S)$ were studied in \cite{T} by using the conjugacy scheme of the group $S_n.$ 
%However, in contrast with the next two subsections, this scheme is neither induced by a graph nor $Q$-polynomial.
For next paragraph, define $$E_i=\{(x,y) \in S_n^2 \mid d_S(x,y)=i\}$$ for all $i\in \{0,1,\dots, n\}.$
In that range of $i,$ write $v_i=|E_i|.$ Note that $v_i=w_{n-i}.$
If $Y\subseteq S_n$ is any set of permutations its {\bf covering radius} $\rho(Y)$ is defined as
$$ \rho(Y)=\max\{ \min\{  d_S(x,y) \mid y \in Y\}\mid x \in S_n\}.$$
%%%%%%%%%%%%%%%%%%%%%
\subsection{Permutation designs}
%%%%%%%%%%%%%%%%%%%%%%%%%%%%%%%%%%%%%%%%%5
If $D$ is any non void subset of $S_n$ we define its {\em frequencies} as
$$\forall i \in [0..n], \, f_i=\frac{|D^2\cap E_i|}{|D|^2}.$$ Thus $f_0=\frac{1}{|D|},$ and $\sum\limits_{i=0}^n f_i=1.$ Note also that if $D=S_n,$ then $f_i=\frac{v_i}{n!}.$
%Consider the random variable $a_D$ defined on $D^2$ with values in $[0..n]$ which affects to an equiprobably chosen 
%$(x,y)\in D^2$ the only $i$ such that $(x,y)\in E_i.$ Thus the frequencies $f_i=Prob(a_D=i).$ 

{\definition The set $D \subseteq X$ is a {\em $t$-design} for some integer $t$ if

$$\sum_{j=0}^n f_jj^i=          \sum_{j=0}^n \frac{v_j}{n!}j^i.           $$
%$$\EE(a_D^i)=\EE(a_X^i)$$ 
for $i=1,\dots,t.$
}
%Let $z_0,z_1,\cdots,z_d$ denote arbitrary pairwise distinct real numbers.

(Note that trivially $\sum\limits_{j=0}^n f_jj^0=1$ so that we do not consider $i=0.$)
Thus, distances in $t$-designs are very regularly distributed. For a $2$-design, for instance, the average and variance of the distance coincide with that of the whole space.

{\bf Remark:} Our notion of design is a special case of designs in  polynomial spaces of \cite{G}.

\subsection{Orthogonal polynomials}

{\definition We define a scalar product on $\R[x]$ attached to $D$  by the relation
$$\langle f,g \rangle_D=\sum_{i=0}^n f_i f(i)g(i). $$
Thus, in the special case of $D=S_n$ we have
$$\langle f,g \rangle_{S_n}=\frac{1}{n!}\sum_{i=0}^n v_i f(i)g(i). $$}

%The following Proposition is given without proof in \cite[p.26]{Sz}.

We require the so-called Charlier polynomials.

Let $$C_k(x)=(-1)^k+\sum\limits_{i=1}^k(-1)^{k-i}{k \choose i} x(x-1)\cdots(x-i+1).$$ An exponential generating function is
given in \cite[(1.12.11)]{KS} as:
$$e^t(1-t)^x=\sum_{n=0}^\infty C_n(x)\frac{t^n}{n!}. $$  

Thus, for concreteness, $C_0(x)=1,\,C_1(x)=x-1,\, C_2(x)=x^2-3x+1.$

%Denote by $F(\sigma, \tau)$ the number of fixed points of $\sigma\tau^{-1}.$ let $d=n-F$ denote the Hamming distance on the permutation images.
The scalar product attached to the space $(S_n,d_S)$ is then 
 $$\langle f,g\rangle_n=\frac{1}{n!} \sum_{k=0}^n w_{n-k} f(k)g(k).$$
 
 It is remarkable that the following orthogonality relation is not found in the classical treatises \cite{KS,Sz} on orthogonal polynomials. See \cite[Lemma 1]{SXS} for a proof.

{\lem \label{Charlier} The reversed Charlier polynomials $\widehat{C_k(x)}=C_k(n-x)$ satisfy the orthogonality relation $$ \langle\widehat{C_r},\widehat{C_s}\rangle_n=r! \delta_{rs},$$
for $r,s\le n/2,$ where $\delta$ denotes the Kronecker symbol.}

%{\definition For a given $D\subseteq X$ the {\em cumulative distribution function} (c.d.f.) is defined as $$F_D(x)=Prob(a_D\le x)=\sum\limits_{i\le x}f_i.$$}
%%%%%%%%%%%%%%%%%%%%%%%%%%%%%%%%%%%%%%%%%%%%%%%%%%%%%%%%%
\section{Structure theorems}
The following result is derived in \cite{CG}, and in a different language in \cite{SXS}.
{\thm If $D\subseteq S_n$ is a $t$-transitive permutation group then it is a $t$-design in $(S_n,d_S).$ 
If $D\subseteq S_n$ is a $t$-design that is a subgroup of $S_n,$ then it is a $t$-transitive permutation group.
}

We require the following characterization of $1$-designs from \cite{SXS}.
{\lem \label{oned} A subset $D\subseteq S_n$ is a $1$-design in $(S_n,d_S)$  iff $\sum\limits_{j=0}^njf_j=n-1.$ In particular,
this condition is satisfied if we have $n$ permutations at pairwise distance $n$ when $f_1=f_2=\dots=f_{n-1}=0,$ and $f_n=\frac{n-1}{n}.$}

%%%%%%%%%%%%%%%%%%%%%%%%%%%%%%%%%%%%%%%%%%%%%%%%%%%%%%%%%
\section{Main result}
We begin with two Lemmas on the zeros of Charlier polynomials.

{\lem \label{Sz} The polynomial $C_k$ has exactly $k$ real zeros.  }
\begin{proof}
 Direct application of Theorem 3.3.1 of \cite{Sz} to the Charlier polynomials which are orthogonal wrt the probability measure of a Poisson law of parameter one \cite[p.34]{Sz}.
\end{proof}

In the following, we will denote by $x(k)$ the largest zero of $C_k.$ This definition makes sense by Lemma \ref{Sz}.

{\lem \label{root} If $z$ is a zero of $C_k$ for $k>1,$ then $z$ cannot be an integer.}

\begin{proof}
 If $z$ is an integral zero of $C_k(x)$ then $z$ divides $C_k(0)=(-1)^k,$ hence $z=1.$ But $C_k(1)=(-1)^k(1-k)$ which is $\neq 0$ for $k>1.$
\end{proof}

Define the {\bf half-strength} of a $t$-design as $s=\lfloor \frac{t+1}{2}\rfloor.$ The next result, which motivates this note, derives an upper bound on the covering radius of a design of given half-strength.

{\thm If $D$ is a design of $S_n$ of half-strength $s>1,$ then $$\rho(D) < n-x(s).$$}

\begin{proof}
 Define $P_s(x)=\widehat{C}_s(x)/(n-x+x(s)).$ Since $P_s$ is a polynomial of degree $<s,$ we have, by orthogonality
 $$\langle \widehat{C}_s, P_s\rangle_{S_n}=\langle 1, \widehat{C}_s P_s \rangle_{S_n}=0.$$
 The degree of $\widehat{C}_s P_s$ is $s+s-1 \le t,$ we have by definition of a $t$-design 
 $$\langle 1, \widehat{C}_s P_s \rangle_D=\langle 1, \widehat{C}_s P_s \rangle_{S_n}=0.$$
 Since, for a given $\sigma \in S_n,$  the translate $\sigma D$ is also a $t$-design, we can write
 
 $$\langle 1, \widehat{C}_s P_s \rangle_{\sigma D}=0.$$
 
 If we assume, looking for a contradiction, that $\rho(D) \ge n-x(s)$ we see that all terms in the above sum being nonnegative, must be zero. Hence all distances of $\sigma$ to $D$ must be roots of 
 $\widehat{C}_s,$ which is impossible for $s>1,$ by Lemma \ref{root}.
\end{proof}

{\ex Computing the roots of $C_k(x)$ using Wolfram online yields

\begin{itemize}
 \item If $t=2$ then $s=1$ and $\rho(D) < n-1,$ hence $\rho(D) \le n-2,$ since $x(1)=1.$
 \item If $t=3$ or $t=4$ then $s=2$ and  $\rho(D) \le n-3,$ since $x(2)\approx 2.616.$
 \item If $t=5$ or $t=6$ then $s=3$ and  $\rho(D) \le n-5,$ since $x(3)\approx 4.115.$
 \item If $t=7$ or $t=8$ then $s=4$ and  $\rho(D) \le n-6,$ since $x(4)\approx 5.544.$
\end{itemize}

In the cases $t=3,5$ our bound coincides with that of \cite{CW}. For $t=2,4,6,7$ it is weaker by one unit.
 }
 
 In general it is known that $x(k)\le k+2\sqrt{k} +1.$ See \cite[Th. 4]{K}.
 
If the strength is small, a direct power moment method is often more effective.

{\thm If $D$ is a $1$-design of $S_n$  then $\rho(D) \le n-1.$}

\begin{proof}
 By Lemma \ref{oned}, we know that the average distance in the shifted design $\sigma D$ is $n-1.$ hence $d(\sigma, D)\le n-1.$ Since $\sigma$ is arbitrary in $S_n,$ the result follows.
\end{proof}

Thus, for $t=1$ also, our bound coincides with that of \cite{CW}.

\section{Conclusion}
In this note, we have studied the covering radius of permutation designs. We have obtained a general upper bound (Theorem 1) on that quantity, dependent on the largest zero $x(t)$ of the Charlier polynomial of degree $t.$
In order to compare Theorem 1 with the bound of \cite{CW}, we would need an asymptotic equivalent of, or a lower bound on $x(t)$ when $t\to \infty.$ We could not find any such result
in the literature of orthogonal polynomials \cite{K,Sz}. This is the main open problem.

{\bf Acknowledgement:} We thank Alexis Bonnecaze and Sam Mattheus for helpful discussions.

%%%%%%%%%%%%%%%%%%%%%%%%%%%%%%%%%%%%%%%%%%%%%%%%%%%%%%

\end{document}